\documentclass[12pt]{article}
\usepackage[T2A]{fontenc}
\usepackage[cp1251]{inputenc}
\usepackage[russian]{babel}
\usepackage{amssymb,amsmath,latexsym,amsthm}
\usepackage{amsbsy,amsfonts,amssymb,mathrsfs,amscd}
\usepackage[mathscr]{eucal}
\usepackage{graphicx}
\usepackage{indentfirst}
\usepackage{verbatim}
\usepackage{enumerate}
\usepackage{amsbib}
\usepackage{cmap}

\theoremstyle{plain}

\newtheorem{theorem}{Теорема}
\newtheorem{corollary}{Следствие}

\newcommand{\RR}{\mathbb{R}} 
\newcommand{\CC}{\mathbb{C}} 
\newcommand{\NN}{\mathbb{N}}

\DeclareMathOperator{\dd}{d\!}
\DeclareMathOperator{\sbh}{sbh}  
\newcommand\pr{\operatorname{pr}} 
\DeclareMathOperator{\clos}{clos} 
\DeclareMathOperator{\Meas}{Meas}
\DeclareMathOperator{\lin}{lin}
\DeclareMathOperator{\loc}{loc}

\begin{document}
{\begin{flushleft} {УДК 517.574 + 517.982.1+517.55 + 517.987.1}\end{flushleft} 

\begin{center}
{\large \bf Нижние огибающие относительно выпуклых подконусов 
субгармонических и плюрисубгармонических функций\footnote{Исследование выполнено за счёт гранта Российского научного фонда  № 22-21-00026, \href{https://rscf.ru/project/22-21-00026/}{https://rscf.ru/project/22-21-00026/}}}
\end{center}

\begin{center}
{\sl Э.~Б.~Меньшикова, Б.~Н.~Хабибуллин}
\end{center}

Часть  функционально-аналитических результатов из \cite{KhaRozKha19} и \cite{Kha21}, где достаточно детально изложена и история вопроса с обширной библиографией, применяется  к двойственному описанию  огибающих из  выпуклых конусов, состоящих из (плюри)субгармонических функций. 

Множества $\NN:=\{1,2, \dots\}$, $\RR$,  $\CC$ \textit{натуральных, вещественных, комплексных чисел,\/} $\NN_0:=\{0\}\cup \NN$ и  $\overline{\RR}:=\RR\cup\{\pm\infty\}$,
где $-\infty:=\inf \RR=\sup \emptyset $ и  $+\infty:=\sup \RR=\inf \emptyset$ для пустого множества $\emptyset$,   рассматриваются с  их естественными структурами.  
$\RR^d$ ---  евклидово пространство над $\RR$  размерности $d\in \NN$
с мерой Лебега $\mathfrak{m}_d$, а  $\bar B_o(r)$  ---  замкнутый шар в $\RR^d$ c центром $o\in \RR^d$ радиуса $r> 0$.  $C(D)$ --- векторное пространство над $\RR$ вещественных непрерывных функций на $D\subset \RR^d$.  Всюду далее $D\subset \RR^d$ --- область,    $ \bar{B}_o(r) \subset  D$, 
$\Meas_0^+(D)$ --- класс всех положительных конечных борелевских мер с компактными носителем в $D$, $\operatorname{sbh}(D)$ --- выпуклый конус всех субгармонических на $D$ функций. Понятия суммируемости интегралов, а также равенств $\overset{\text{п.в.}}{=}$
и неравенств $\leqslant^{\text{п.в.}}$ почти всюду (п.в.)
относятся ниже именно к мере Лебега   $\mathfrak{m}_d$.
 \textit{Постоянная} $c\in \overline \RR$ рассматривается и как \textit{функция, тождественно равная\/ $c$.}     Следующий результат решает
\cite[п.~2.3, задачи 1--2]{KhaRozKha19}, \cite[п.~1.2.3, задачи 1--2]{Kha21} для широкого круга выпуклых подконусов в $\sbh(D)$.  

\begin{theorem}\label{th1}
Пусть  выпуклый подконус $H\subset \sbh (D)$ содержит  постоянные\/   $0$ и\/  $-1$, и    для любой  локально ограниченной сверху на $D$ последовательности  $(h_k)_{k\in \NN}$ функций $ h_k\in H$ полунепрерывная сверху  регуляризация   
поточечного верхнего  предела этой  последовательности  принадлежит $H$.    
Тогда для любой определённой  п.в. функции $f$ на $D$, равной п.в. некоторой функции из   $C(D)+ H-H$, выполнено равенство
\begin{align}\label{supinf}
\sup&\biggl\{\int_{\bar{B}_o(r)}
h\dd {\mathfrak m}_d\Bigm| -\infty  \neq h\in H,   h{\leqslant^\text{п.в.}}  f \text{ на $D$}\biggr\}=\inf_{\mu \in J_H} \int_{D}f\dd \mu,
\\
\intertext{где } J_H&:=\biggl\{ \mu\in \Meas_0^+(D)\biggm|  \int_{ \bar{B}_o(r)} 
h\dd {\mathfrak m}_d\leqslant \int_{D}h\dd \mu
\text{ для всех $h\in H$} \biggr\}. 
\label{J}
\end{align}
В частности,  существование функции $h\in H$ с требованиями 
$$-\infty \neq h{\leqslant^{\text{п.в.}}} f\quad\text{на $D$}
$$
 равносильно тому, что 
$\inf$ справа в \eqref{supinf} по классу всех мер $\mu\in J_H$ не равен $-\infty$.
\end{theorem} 
Примерами конусов $H$, удовлетворяющих  условиям теоремы  \ref{th1},
являются конус $\operatorname{sbh}(D)$; 
конус $\operatorname{psbh}(D)\subset \sbh(D)$ всех\textit{ плюрисубгармонических  на\/} $D\subset \CC^d$ функций, если  $\CC^d$ отождествить с $\RR^{2d}$; выпуклые подконусы этих конусов, выделяемые условиями отрицательности $\leqslant 0$ на фиксированном открытом подмножестве в $D$, или инвариантные относительно определённых
преобразований $D$ и мн. др.   Наиболее важна в теореме \ref{th1} широкая 
возможность выбора $f$ из   $H-H+C(D)$. Случай $f\in C(D)$ ранее 
был полностью разобран в  \cite[следствие 8.1]{KhaRozKha19},  \cite[следствие 3.2.1]{Kha21}, 
\cite[теорема 7.2]{Kha01}.   При $H:=\operatorname{psbh}(D)$
\textit{для равенств вида} \eqref{supinf} в  
\cite{Pol91}--\cite{Kuz18}
от $f$ всегда  требовалась локальная ограниченность сверху на области $D$. Но  при  крайне актуальном для дальнейших применений варианте $f\in H-H$ локальная ограниченность сверху скорее  не  выполняется, поскольку функции из $-H$ могут быть не ограничены сверху даже на всюду плотном в  $D$ подмножестве. То же самое в отношении $H$ по локальной ограниченности снизу. 

\begin{corollary}\label{cor1} 
Пусть $D\subset \CC^d$ --- область, $v$ и  $ M\neq -\infty$ --- пара плюрисубгармонических функций на $D\subset \mathbb{C}^d$. Существование    $h\neq -\infty$ из $\operatorname{psbh}(D)$ c ограничением 
 $$
v+h\leqslant M\quad\text{на $D$}
$$
 равносильно  условию 
$$
\sup\limits_{\mu\in J_{\operatorname{psbh}(D)}}\int (v-M)\dd \mu<+\infty,
$$
где супремум в последнем неравенстве берётся  по классу $J_H$, определённому в \eqref{J},  при $H:=\operatorname{psbh}(D)$.
 \end{corollary}

\begin{proof}[Доказательство теоремы\/ {\rm \ref{th1}}] Нам потребуются следующие понятия, а также   не\-к\-о\-т\-о\-рые факты  о них  из \cite{KhaRozKha19}--\cite{Kha01}.  

Упорядоченное векторное пространство $(X, \leq)$ над $\RR$ c отношением порядка $\leq$ 
--- \textit{векторная решётка,\/} если существует  \textit{точная верхняя грань\/} 
$X\text{-}\sup F\in X$ для любого конечного $F\subset X$. 
Для пары векторных решёток $X$ и $Y$ через $\lin^+ Y^X$  обозначаем \textit{выпуклый  конус  линейных положительных,\/} или \textit{возрастающих, функций\/} $l\colon X\to Y$.
Пусть $(X_n)_{n\in \NN_0}$ --- последовательность  \textit{векторных решёток\/}  $(X_n, \le_n)$, 
 $$\prod X_n:=\prod_{n=0}^{\infty} X_n$$ --- их произведение, и для  $$x=(x_n)_{n\in \NN_0}\in \prod X_n$$ 
через 
$$
\pr_n x=x_n\in X_n,
$$
обозначаем проекцию $x$ на 
Имеет место неравенство 
$$x\leqslant x'\quad\text{в}\quad \prod X_n, 
\quad\text{если   $\pr_n x\leq_n \pr_n x'$ для каждого  $n\in{\NN_0}$.}
$$
Пусть 
$$p_n\in \lin^+ X_n^{X_{n+1}}$$  
и   
$$
X_n\text{-}\sup p_n(F_{n+1})=p_n\bigl(X_{n+1}\text{-}\sup F_{n+1}\bigr)\quad 
$$ 
{\it для  каждого конечного\/ $F_{n+1}\subset X_{n+1}$.}
Тогда 
$$
\projlim X_np_n :=\bigl\{x\in \prod X_n\bigm|  \pr_n x = p_n (\pr_{n+1} x)\quad \text{\it при всех  $n\in \NN_0$} \bigr\}
$$ с тем же отношением порядка  $\leqslant$ ---  векторная решётка,
называемая   {\it проективным пределом последовательности  $(X_n)_{n\in \NN_0}$ по  $(p_n)_{n\in \NN_0}$.\/} Можно считать \cite[предложение 3.1]{KhaRozKha19}, \cite[предложение 2.1.1]{Kha21},  что  
$$
\pr_n X:=\bigl\{\pr_nx\bigm| x\in X\bigr\}=X_n\quad \text{для любого $n\in \NN$}.
$$ 
Подмножество $B\subset X$ \textit{ограничено снизу (сверху) в\/} $X$,  если 
существует $x\in X$, для  которого $x\leqslant b$ (соответственно $b\leqslant x$) для всех $b\in B$, и $B$ \textit{ограничено в\/} $X$, если $B$ ограничено  и  снизу, и сверху в $X$. 

\begin{theorem}[{\cite[теорема 2, следствия 6.1 и 3.1]{KhaRozKha19}, \cite[теорема 2.4.1, следствия 2.4.1 и 2.1.1]{Kha21}}]\label{th2}
 Пусть $H_*\subset X:=\projlim  X_np_n$ --- выпуклый конус и\/ $0\in H_*$, 
а для любой ограниченной в $X$ последовательности  $(h^{(k)})_{k\in \NN}$, $h^{(k)}\in H_*$,
 существует принадлежащий $H_*$ верхний предел 
\begin{equation}\label{limh}
\limsup\limits_{k\to \infty} h^{(k)} :=\inf\limits_{n\in \NN}\sup\limits_{k\geqslant n} h^{(k)} \in H_*.
\end{equation}

Пусть $S\subset X$ --- векторное подпространство, содержащее  $H_*$, и при  каждом $n\in \NN_0$ для любого $s_n\in \pr_n S$ найдётся  $h_n\in \pr_n H_*$, удовлетворяющее неравенству $h_n\leq_n s_n$ в $X_n$. 

Пусть для   $q_0\in  \lin^+\RR^{X_0}$  и для  суперпозиции  $$q:=q_0\circ \pr_0\in \lin^+\RR^{X}$$
 при любой  убывающей  в $X$ последовательности $(h^{(k)})_{k\in \NN}$, $h^{(k)}\in H_*$,   при конечности $$\inf\limits_{k\in \NN} q(h^{(k)})\in \RR$$ эта  последовательность  $(h^{(k)})_{k\in \NN}$ ограничена снизу в $X$ и  $$q\bigl(\inf_{k\in \NN} h^{(k)}\bigr)\geqslant \inf_{k\in \NN} q(h^{(k)}).$$ 

Тогда  для каждого  $s\in S$ величина $$\sup \bigl\{ q(h) \bigm| H_*\ni h\leqslant s\bigr\}\in \overline \RR$$
 равна величине  
\begin{equation}\label{eq5_1}
\hspace{-2mm}
\inf \biggl\{(l_n\circ \pr_n)(s) \Bigm| n\in \NN_0,  l_n\in \lin^+ \RR^{\pr_nS},     
q(h)\underset{\forall h\in H_*}{\leqslant} (l_n\circ \pr_n)(h) \biggr\}\in \overline \RR.
\end{equation}
\end{theorem}
Для $D\subset \RR^d$  выберем   исчерпание последовательностью\/ $(D_n)_{n\in \NN_0}$ областей   $D_n\subset \RR^d$, для которого $\overline B_o(r)\subset D_0$, \textit{замыкание\/} $\clos D_n$ области $D_n$ содержится в области $D_{n+1}$ при каждом $n\in \NN_0$ и $D=\bigcup_{n\in \NN_0} D_n$. 

Для $n\in \NN_0$ рассмотрим пространство $X_n:=L^1(\clos D_n)$ суммир\-у\-е\-мых на $\clos D_n$ функций с отношением поточечного предпорядка ${\leq}_n^{\text{п.в.}}$, факторизацию  которого  по отношению $\overset{\text{п.в.}}{=}$ обозначим  через $X_n$, где  ${\leq}_n^{\text{п.в.}}$  
уже отношение порядка.  В качестве линейных положительных функций $p_n\in \lin^+ X_n^{X_{n+1}}$ выберем \textit{сужения <<функций>> из $X_{n+1}$ на $\clos D_{n+1}$,\/} которые становятся уже векторами  из  
 $X_{n}$.  Проективный предел векторных решёток $\projlim X_np_n$ 
в этом случае  --- это факторизованное   по отношению  $\overset{\text{п.в.}}{=}$  пространство  локально суммируемых на $D$ функций,  с  отношением порядка 
$\leqslant^{\text{п.в.}}$, которое обозначаем через   $L_{\loc}^1(D)$. 

Для выпуклого конуса $H\subset \sbh(D)$ положим $$H_*:=H\setminus \{-\infty\}
\subset \sbh(D)\setminus \{-\infty\}\subset L_{\loc}^1(D).$$   
Полунепрерывная сверху регуляризация верхнего предела последоват\-е\-ль\-н\-ости субгармонических функций на области, --- если этот верхний предел  не равен $-\infty$, ---  
c одной стороны даёт субгармоническую функцию, а с другой 
отличается от верхнего предела разве что на   множестве  нулевой $\mathfrak{m}_d$-меры, и даже меньшем.  Поэтому для этого конуса $H_*$ выполнено условие, завершающееся  соотношением \eqref{limh}.  

Положим 
$$
S:=C(D)+H_*-H_*\subset L_{\loc}^1(D).$$
 Очевидно, $H_*\subset S$. 
Пусть $s_n\in S$, т.е. 
$$
s_n=g_n+h_n-h_n',\text{ где $g_n\in C(\clos D_n)$,}
$$
 а $h_n\in \pr_n H_*$ и 
$h_n'\in \pr_nH_*$  --- сужения  на $\clos D_n$ функций из $H_*$.
Тогда существуют  \textit{положительные   числа\/} $c$ и $c'$,
 для которых  $g_n\geqslant -c$ на $\clos D_n$ и $h_n'\leq_n c'$ на $\clos D_n$. 
Следовательно, $$
h_n-c-c'\leq_n s_n,\text{ где $h_n\in \pr_n H_*$,}
$$ а также 
отрицательная постоянная   $-c-c'=(c+c')(-1)$ принадлежит $\pr_nH$, поскольку по  условию 
$-1\in H$.  Таким образом, выполнены условия теоремы, касающиеся  подпространства $S$.

В качестве $q_0$ в теореме \ref{th2} рассмотрим сужение меры $\mathfrak{m}_d$ на $\overline B_o(r)$ в том смысле, что 
\begin{equation}\label{q0f}
q_0(f_0):=\int_{\overline B_o(r)} f_0\dd \mathfrak{m}_d\in \RR
\quad\text{для всех $f_0\in X_0=\pr_0 L_{\loc}^1(D)$}.
\end{equation}
Функция $u\colon D\to \overline \RR$ \textit{почти субгармоническая\/} на $D$, если она   п.в.  совпадает с субгармонической функцией. Для \textit{произвольной убывающей п.в. последовательности\/} $(h^{(k)})_{k\in \NN}$ почти субгармонических на $D$ функций $h^{(k)}$  условие 
$$
\inf\limits_{k\in \NN} q(h^{(k)})\in \RR\quad\text{для $q:=q_0\circ \pr_0$}
$$
 означает, что  
\begin{equation}\label{infi}
\inf\limits_{k\in \NN}\int_{\overline B_o(r)}h^{(k)}\dd  \mathfrak{m}_d>-\infty.
\end{equation}
Отсюда предел этой последовательности даёт почти субгармоническую функцию на $D$.
Тем более это  верно для убывающей последовательности $(h^{(k)})_{k\in \NN}$  из конуса $H_*\subset \sbh(D)\setminus \{-\infty\}$. 
Верхний предел \eqref{limh} убывающей последовательности --- это точная нижняя грань этой последовательности. Поэтому по условию теоремы \ref{th1} для последовательностей 
 $(h^{(k)})_{k\in \NN}$ при условии   \eqref{infi}  получаем  $$-\infty\neq\inf_{k\in \NN} h^{(k)}\in H.$$
 В частности,  при \eqref{infi}  убывающая последовательность $(h^{(k)})_{k\in \NN}$ из $H_*$, очевидно, ограниченная сверху функцией $h^{(1)}$,  ограничена и снизу функцией $$\inf_{k\in \NN} h^{(k)}\in H_*.$$
  При этом можем считать все 
функции $h^{(k)}$ полунепрерывными сверху. Для убывающей последовательности таких функций $\inf\limits_{k\in\NN}$  можно внести под знак интеграла:
\begin{equation*}
\inf\limits_{k\in \NN}\int_{\overline B_o(r)}h^{(k)}\dd  \mathfrak{m}_d
=\int_{\overline B_o(r)} \inf\limits_{k\in \NN} h^{(k)}\dd  \mathfrak{m}_d.
\end{equation*}
Это означает, что выполнено 
$$
q\bigl(\inf_{k\in \NN} h^{(k)}\bigr)\geqslant \inf_{k\in \NN} q(h^{(k)})\quad\text{для $q:=q_0\circ \pr_0$.}
$$ В итоге из условий теоремы \ref{th1}
следуют  все условия теоремы \ref{th2}. Легко видеть, 
что $$\sup \bigl\{ q(h) \bigm| H_*\ni h\leqslant s\bigr\}\in \overline \RR$$
 --- 
это левая  часть в \eqref{supinf}. Убедимся, что \eqref{eq5_1} в рассматриваемой ситуации     --- это правая часть в \eqref{supinf}. 

Пусть, как в \eqref{eq5_1},  $l_n\in \lin^+ \RR^{\pr_nS}$, где $S=C(D)+H_*-H_*$.
Тогда  $$l_n\in \lin^+ \RR^{\pr_n C(D)}=\lin^+ \RR^{C(\clos D_n)},
$$
 откуда 
по теореме Рисса $l_n$ на $C(\clos D_n)$ реализуется как некоторая 
положительная конечная мера Бореля $\mu$ на $D$ c компактным носителем в $\clos D_n$. 
Требование из  \eqref{eq5_1} вида 
 $$q(h)\leqslant (l_n\circ \pr_n)(h)\quad\text{при всех  $h\in H_*$}
$$
для $q=q_0\circ \pr_0$ согласно \eqref{q0f} в терминах меры $\mu$ можно записать как
$$
\int_D h\dd \mathfrak{m}_d\leqslant \int_D h\dd \mu
\quad\text{при  всех $h\in H_*$}. 
$$
Последнее влечёт за собой конечность  интегралов   
$$
\int_D h\dd \mu\in \RR
\quad\text{для  всех $h\in H_*$.}
$$ 
Следовательно, полученные таким образом 
 меры $\mu\in \Meas^+_0(D)$  корректно определены на 
$S=C(D)+H_*-H_*$ и пробегают в точности $J_{H_*}=J_{H}$ из \eqref{J}, 
поскольку удаление  постоянной $-\infty$ из  $H$ в  \eqref{supinf}--\eqref{J} ничего не меняет. 
Таким образом равенство \eqref{supinf} доказано, а заключительное утверждение теоремы \ref{th1} --- очевидное следствие этого равенства. 
\end{proof}

\begin{proof}[Доказательство следствия\/ {\rm \ref{cor1}}] 
Достаточно рассмотреть $v\neq -\infty$. 
Для любой пары функций $u_1,u_2\in \sbh(D)\setminus \{-\infty\}$ из неравенства 
$u_1{\leqslant^{\text{п.в.}}}u_2$ на $D$ следует 
$u_1\leqslant u_2$ всюду на $D$. Поэтому  существование функции 
 $h\in \operatorname{psbh}(D)\setminus \{-\infty\}$ c ограничением 
 $$
v+h\leqslant M \quad\text{на $D$}
$$
  равносильно существованию  
 функции   $h\in \operatorname{psbh}(D)\setminus \{-\infty\}$ c ограничением 
 $$
h{\leqslant^{\text{п.в.}}} (M-v)=:f\quad\text{на $D$,}
$$
 а это эквивалентно конечности левой части 
\eqref{supinf} при $H:= \operatorname{psbh}(D)$. Последнее по заключительному утверждению 
теоремы  \ref{th1} равносильно условию 
$$
\inf_{\mu \in J_H} \int_{D}f\dd \mu>-\infty,
$$ 
что эквивалентно  условию 
$$
\sup\limits_{\mu \in J_H} \int_{D}(v-M)\dd \mu<+\infty.
$$  
Следствие доказано.\end{proof}

\vskip 3mm

Институт математики  с вычислительным центром  Уфимского федерального исследовательского центра РАН

E-mail: {\tt enzha@list.ru  \quad khabib-bulat@mail.ru}
\end{document}